\documentstyle[12pt]{article}
\begin{document}

\title{ \Large {Grassmann Manifold $G(2,8)$  and  Complex Structure on $S^6$}
 \footnotetext{
\mbox{}\quad  \ \ {\it Subject classification}. \ 14M15, 53C15,
53C27. }}
\author{Jianwei Zhou}
\date{\small Department of Mathematics, Suzhou University, Suzhou
215006, P.R. China }

\maketitle
\begin{abstract} In this paper, we use  Clifford algebra and the spinor calculus to
 study the complex structures on Euclidean
space $R^8$ and the spheres $S^4,S^6$. By the spin representation of $G(2,8)\subset Spin(8)$ we show that the
Grassmann manifold $G(2,8)$ can be looked as the set of orthogonal complex structures on $R^8$. In this way, we
show that $G(2,8)$ and $CP^{3}$ can be looked as twistor spaces of $S^6$ and $S^4$ respectively. Then we show that
there is no almost complex structure on sphere $S^4$ and there is no orthogonal complex structure on the sphere
$S^6$.

\end{abstract}

\baselineskip 15pt
\parskip 5pt

\vskip 1cm \centerline{\bf \S 1. Introduction }

\vskip 0.3cm

In this  paper, we study the complex structures on Euclidean space
$R^8$ and the spheres $S^4,S^6$ respectively. In the study we use
 Clifford algebra and the spinor calculus. The main references
are [4], [5], [6] and [7].

 Let $\bar e_1, \bar e_2, \cdots, \bar e_8$ be a fixed
orthonormal basis of $R^8$, the Clifford product on Clifford
algebra $C\ell_8$ be determined by the relations: $$\bar e_B\bar
e_C+\bar e_C\bar e_B =-2\delta_{BC}, \ \ B,  C=1,  2, \cdots, 8.$$
As in [6], let $A_8=Re[(\bar e_1+\sqrt{-1}\bar e_2)\cdots(\bar
e_7+\sqrt{-1}\bar e_8)], \ \beta_8=\bar e_1\bar e_3\bar e_5\bar
e_7, \ A=A_8(1+\beta_8)$. The space $V=C\ell_8\cdot A=V^+\oplus
V^-$ is an irreducible module over $C\ell_8$ and the spinor spaces
$V^+=C\ell_8^{even}A, \ V^-=C\ell_8^{odd}A$ are generated by
$\alpha_B= \bar e_1\bar e_BA$ and $\alpha_{B+8}=\bar e_B A$
respectively, $B=1, \cdots, 8$. For any $x\in C\ell_8$,
$$x(\alpha_1,\cdots,\alpha_8,\alpha_9,\cdots,\alpha_{16})=
(\alpha_1,\cdots,\alpha_8,\alpha_9,\cdots,\alpha_{16})\Phi(x)$$
defines an algebra isomorphism $\Phi\colon\; C\ell_8\to R(16)$.

As shown in [6], for $x\in C\ell_8$, we have
$\Phi(\alpha(x^t))=(\Phi(x))^t$. Then for any $g\in Spin(8), \
\Phi(g)=\left(
\begin{array}{cc} B & {} \\ {} & C
\end{array}
 \right),$ we have $B,C\in SO(8).$ If $g^2=-1, \ B$ and $C$ define
orthogonal almost complex structures on $R^8$.

In this paper, we use  isomorphism $\Phi\colon\; C\ell_8\to R(16)$
to construct homeomorphism  $\Phi^*$ from Grassmann manifold
$G(2,8)$ to the set of oriented orthogonal complex structures on
$R^8$. Furthermore, there are two fibre bundles $\tau \colon\;
G(2,8)\to S^{6}$ and $\tau_1\colon\; CP^{3}\to S^{4}$ defined
naturally. We show that restricting the homeomorphism $\Phi^*$ on
the fibres of these fibre bundles respectively, we get the sets of
complex structures on the tangent spaces of $S^6$ and $S^4$
respectively. By definition on [4] p.339, $G(2,8)$ and $CP^3$ are
the twistor spaces on $S^6$ and $S^4$ respectively.

Then the almost complex structure on the sphere $S^6$ and $S^4$
are determined by  sections $f\colon\; S^6\to G(2,8)$ and $f: \
S^4\to CP^3$ respectively. By the cohomology groups  of  $CP^3$
and $S^4$, we show that there in no almost complex on $S^4$.

 $G(2,8)$ is a Kaehler manifold, in
\S 3 we show that the almost complex structure on $S^6$ is integrable if and only if the map $f\colon\; S^6\to
G(2,8)$ is holomorphic. Then $S^6$ is also a Kaehler manifold if $f\colon\; S^6\to G(2,8)$ holomorphic. These
shows there is no orthogonal complex structure on the sphere $S^6$.

For the complex structures on $R^8,S^6, S^4$ see also [1] and [3],
p.159, p. 281.

\vskip 1cm \centerline{\bf \S 2. The complex structures on $R^8$}

\vskip 0.3cm

$G(2, 8)$ is the Grassmann manifold formed by all oriented
$2$-dimensional subspaces of $R^8$, any $x\in G(2,8)$ can be
represented by $e_1e_2$ or $e_1\wedge e_2$, where $e_1, e_2$ is an
oriented orthonormal basis of $x$. Then $G(2, 8)$ can be viewed as
a subspace of $Spin(8)$ or $\bigwedge^2(R^8)$.  Let
$$M=\{ A\in SO(8) \ |
\ A^2=-I\} \approx SO(8)/U(4)$$ be the set of oriented orthogonal
complex structures on $R^8$.

{\bf Theorem 2.1} \ The map $\Phi^*\colon\; G(2,8)\to M$ defined
by
$$x(\alpha_1,\cdots,\alpha_8)=
(\alpha_1,\cdots,\alpha_8)\Phi^*(x), \ \ x\in G(2,8), $$  is a
diffeomorphism. $\Phi^*(x)$ is the spin representation of $x$.

 {\bf Proof} \  For any $x\in G(2,8)$,
\ $\Phi(x)=\left(
\begin{array}{cc} B & {} \\ {} & C
\end{array}
 \right),$
where $B,C\in SO(8)$. From $x\cdot x=-1$, we have $BB=-I, \ B$
defines  a complex structure on $R^8\cong V^+$. By Proposition 4.3
of [7],  the map $x\in G(2,8)\mapsto B$ is a monomorphism. The
proposition follows from
$$G(2,8)=\{A(\bar e_1)A(\bar e_2) \ | \ A\in SO(8)\}= \{g(\bar
e_1\bar e_2)g^t \ | \ g \in Spin(8)\},$$
$$\Phi(g(\bar
e_1\bar e_2)g^t)=\Phi(g)\Phi(\bar e_1\bar e_2)\Phi(g)^t,$$ and
$$M=\{ABA^t \ | \ A\in SO(8)\},$$
where $B\in M$ is a fixed element.

 For any $v\in R^8$
there is $P_v\in SO(8)$ such that $\Phi(v) = \left(
\begin{array}{cc} {} & P_v \\ -P_v^t & {}
\end{array}
\right)$. Then for any $x=e_1e_2\in G(2,8), \
\Phi(x)=-P_{e_1}P_{e_2}^t$. Furthermore,  $P_v\in M$ if $v\in
S^7,v\perp \bar e_1$.

Denote $J_x\colon\; R^8\to R^8$ the complex structure defined by
$x\in G(2,8)$. For any $v=\sum\limits_{i=1}^8 \ v^i\bar e_i\in
R^8, \ \bar e_1vA = \sum\limits_{i=1}^8 \ v^i\alpha_i$, we have
$$\bar e_1 J_xv A = x\bar e_1vA.$$

It is easy to compute
$$\Phi^*(\bar e_1\bar e_2)=P_{\bar e_2}=\left( \begin{array}{cccccccc}
0  & 1  & &  &  &  &  &  \\
-1  &0  &  &  &  &  &  &  \\
  &  &0  & -1  & &  &  &  \\
 & & 1  &0  &      &   &   \\
& &  &   &0  & -1  &  &   \\
   &   &   &   & 1  &0  &   &   \\
   &  &   &   &   &   &0  & -1  \\
  &  &   &   &  &   & 1  &0  \\
\end{array}\right).$$
The complex structure $J_{x}=\Phi^*(x)$ acts on the left of $R^8$.

As shown in [5] or [7], for any $x\in G(2,8)$ there is $v\in
S^6=\{v\in S^7 \ | \ v\perp \bar e_1\}$ such that $ xA=\bar e_1vA,
\ x\mapsto v$ defines  a  fibre bundle $\tau :  G(2, 8) \to S^6$.
For any $v\in S^6$, we have $\tau ^{-1}(v)=\{uJ_vu \ | \ u\in S^7
\}$, where $J_v$ is a complex structure on $R^8$, see [7] or \S 3
below. The fibres of $\tau $ are diffeomorphic to the complex
projective space $CP^3$.

{\bf Theorem 2.2} \ For any $v\in S^6, \ \Phi^*(\tau ^{-1}(v))$ is
the set of orthogonal almost complex structures on the tangent
space $T_vS^6$. Then $G(2,8)$ is the twistor space of $S^6$ and
any almost complex structure on $S^6$ is defined by a section $f:
\ S^6\to G(2,8)$.

{\bf Proof} \ For any $v\in S^6, \ T_vS^6=\{ X\in R^8 \ | \ X\perp
\bar e_1, v\}$. On the other hand, $x\in \tau ^{-1}(v)$ if and
only if
$$xA=-x\alpha_1=\bar e_1vA, \ \ x\bar e_1vA=xxA=-A.$$
Then the subspace of $R^8$ generated by $\bar e_1,v$ is invariant
under the map  $J_x = \Phi^*(x)$ for any $x\in \tau ^{-1}(v)$ and
$J_x$ gives a complex structure on $T_vS^6$. By Theorem 2.1, all
complex structures on the tangent space  $T_vS^6$ can be obtained
in this way.

By this theorem,  any $x\in \tau ^{-1}(v)$ gives a complex
structure  $J_{x}\colon\; T_vS^6\to T_vS^6$. For any $X\in T_vS^6,
\ Y=J_{x}X$ is defined by $x\bar e_1 XA=\bar e_1YA.$

By definition on [4] p.339, $G(2,8)$ is a twistor space of $S^6$,
 any almost complex structure on $S^6$ is defined by a section
$f\colon\; S^6\to G(2,8)$. In \S 3, we shall make further study
and show that there is no integrable complex structure on $S^6$.

As is well-known, there is a map $CP^3\to S^4$. This map  can also
be constructed by the Clifford algebra $C\ell_8$. In the following
we consider $\tau ^{-1}(\bar e_3)\approx CP^3 $ as an example.

Let $J$ be a complex structure on $R^8$ defined by $J\bar
e_{2i-1}=\bar e_{2i}, \ J\bar e_{2i}=-\bar e_{2i-1}, \ i=1,2,3,4$.
By Lemma 3.6, Proposition 3.7 of [7], for any $x \in
\pi^{-1}(\bar{e}_3)$, there is a vector $v$ such that
$$xA_8\beta_8 = \bar{e}_1vA_8 + \bar{e}_1(\bar{e}_3 -
v)A_8\beta_8,$$ and $v \perp \bar{e}_1, \bar{e}_2, \bar{e}_4$, \
$|\bar{e}_3 - 2v | = 1$. From this equation we have
$$xA_8(1+\beta_8)= \bar e_1\bar e_3A_8(1+\beta_8),$$
$$xA_8(1-\beta_8)= \bar e_1(\bar e_3-2v)A_8(1-\beta_8).$$

Let $S^4=\{t=\bar{e}_4 - 2Jv\in S^7 \ | \ v \perp \bar{e}_1,
\bar{e}_2, \bar{e}_4\}$ be a unit sphere in $R^8$. The maps
$x\mapsto\bar e_3-2v\mapsto\bar e_4-2Jv$ define a map
$$\tau_1\colon\; \tau ^{-1}(\bar
e_3)\approx CP^3\to S^4.$$ By Proposition 2.6 of [7], the 2-form
of  $xA_8\beta_8 = \bar{e}_1vA_8 + \bar{e}_1(\bar{e}_3 -
v)A_8\beta_8$ is a calibration and its contact set diffeomorphic
to $CP^1$ which is also the set $\tau_1^{-1}(t)$. Then $\tau_1: \
\tau ^{-1}(\bar e_3)\to S^4$ is a fibre bundle.

 {\bf Theorem 2.3} \  Restricting the
map $\Phi^*\colon\; G(2,8)\to M$ on the fibre $\tau_1^{-1}(t)$ of
$\tau_1$, we have all complex structures on the tangent space
$T_tS^4$. Then $CP^3$ is the twistor space of $S^4$.

{\bf Proof}  \ It is easy to see  $A_8(1-\beta_8)\bar e_1\bar
e_2=-\bar e_1\bar e_2A_8(1+\beta_8)=-\alpha_2$,  from
$xA_8(1-\beta_8)= \bar e_1(\bar e_3-2v)A_8(1-\beta_8)$ and $\bar
e_1\bar e_2A_8(1+\beta_8)=vJvA_8(1+\beta_8)$ for any $v\in S^7$,
we have
\begin{eqnarray*} x\alpha_2 & = & \bar e_1(\bar e_3-2v)\bar e_1\bar e_2A_8(1+\beta_8)\\
& =& \bar e_1(\bar e_3-2v)\bar e_3\bar e_4A_8(1+\beta_8) \\
& = & -\bar e_1\bar e_4A_8(1+\beta_8)- 2  \bar e_1v
vJvA_8(1+\beta_8) \\
&=& - \bar{e}_1 (\bar{e}_4 - 2Jv) A_8(1+\beta_8).
\end{eqnarray*}

Thus for any $x\in G(2,8), \ xA_8\beta_8 = \bar{e}_1vA_8 +
\bar{e}_1(\bar{e}_3 - v)A_8\beta_8$ if and only if the following
identities hold.
$$x \alpha_1 =- x A_8(1+\beta_8) =   -\alpha_3, $$
$$x \alpha_3  =\alpha_1,$$
$$x\alpha_2 = - \bar{e}_1 (\bar{e}_4 - 2Jv) A_8(1+\beta_8),$$
$$x \bar{e}_1 (\bar{e}_4 - 2Jv) A_8(1+\beta_8)=\alpha_2.$$
These shows that the subspace of $R^8$ generated by $\bar e_1,\bar
e_2, \bar e_3, t=\bar{e}_4 - 2Jv$ is invariant under the action
defined by $J_x$. Thus $J_x$ defines a complex structure on
$T_tS^4= \{ \ w \in R^8 \ | \ w\perp \bar{e}_1, \bar{e}_2,
\bar{e}_3,  t \}$.

Then any almost complex structure on $S^4$ defines a section of
the fibre bundle $\tau_1\colon\;  \tau ^{-1}(\bar e_3)\to S^4$. We
have proved Theorem.

In the following we study the fibres of $\tau_1\colon\; \tau
^{-1}(\bar e_3)\approx CP^3\to S^4$.

First assuming $t = \cos \theta \bar{e}_4+ \sin \theta \bar{e}_6
\in S^4$. From $\bar{e}_3 - 2v = \cos \theta \bar{e}_3+ \sin
\theta \bar{e}_5 $, we have
$$v = \sin \frac{\theta}{2}(\sin \frac{\theta}{2}e_3 - \cos
\frac{\theta}{2}e_5),\ \ e_3 - v = \cos \frac{\theta}{2}(\sin
\frac{\theta}{2}\bar{e}_3 + \sin \frac{\theta}{2}\bar{e}_5).$$
Then
\begin{eqnarray*}&& \bar{e}_1vA_8 + \bar{e}_1(\bar{e}_3 - v)A_8\beta_8 \\
& & =\bar{e}_1\sin \frac{\theta}{2}(\sin \frac{\theta}{2}e_3 -
\cos \frac{\theta}{2}e_5)\bar e_1\bar e_3\bar e_5\bar
e_7A_8\beta_8 \\
&& \quad + \bar{e}_1 \cos \frac{\theta}{2}(\sin
\frac{\theta}{2}\bar{e}_3 + \sin
\frac{\theta}{2}\bar{e}_5)A_8\beta_8 \\
 & & = (
\cos \frac{\theta}{2}\bar{e}_1 + \sin \frac{\theta}{2}\bar{e}_7)
(\cos \frac{\theta}{2}\bar{e}_3 + \sin
\frac{\theta}{2}\bar{e}_5)A_8\beta_8.
   \end{eqnarray*}
By Proposition 2.6 of [7], the calibration defined by
$\bar{e}_1vA_8 + \bar{e}_1(\bar{e}_3 - v)A_8\beta_8$ is
$$( \cos \frac{\theta}{2}\bar{e}_1 + \sin \frac{\theta}{2}\bar{e}_7)
(\cos \frac{\theta}{2}\bar{e}_3 + \sin \frac{\theta}{2}\bar{e}_5)-
( \cos \frac{\theta}{2}\bar{e}_2 + \sin \frac{\theta}{2}\bar{e}_8)
(\cos\frac{\theta}{2}\bar{e}_4 + \sin
\frac{\theta}{2}\bar{e}_6).$$ This is a special Lagrangian
calibration and the  set $\tau_1^{-1}(t)$ is the contact set of
this calibration which  diffeomorphic to $CP^1$. For any
$x\in\tau_1^{-1}(t)$, we have
  \begin{eqnarray*} xA_8\beta_8 & = & \bar{e}_1vA_8 + \bar{e}_1(\bar{e}_3 - v)A_8\beta_8 \\
  &=&
(\cos \frac{\theta}{2}\bar{e}_1 + \sin \frac{\theta}{2}\bar{e}_7)
(\cos \frac{\theta}{2}\bar{e}_3 + \sin
\frac{\theta}{2}\bar{e}_5)A_8\beta_8.  \end{eqnarray*}

Lev $V$ be a subspace of $R^8$ generated by
$$\cos \frac{\theta}{2}\bar{e}_1 + \sin \frac{\theta}{2}\bar{e}_7,
\ \ \cos \frac{\theta}{2}\bar{e}_3 + \sin
\frac{\theta}{2}\bar{e}_5, \ \ \cos \frac{\theta}{2}\bar{e}_2 +
\sin \frac{\theta}{2}\bar{e}_8, \ \ \cos\frac{\theta}{2}\bar{e}_4
+ \sin \frac{\theta}{2}\bar{e}_6.$$ Then any element of
 $\tau_1^{-1}(t)$ can be represented by $vJ_{\bar e_3}v, \ v\in
 V$.

For general $u\in S^4$, we can choose an element $G\in U(2)$ such
that $G(u)=t= \cos \theta \bar{e}_4+ \sin \theta \bar{e}_6,$ where
$U(2)\subset Spin_7$ is the unitary group which fixed the elements
$\bar{e}_1, \cdots, \bar{e}_4$. Since $A_8(1+\beta_8)$ is
invariant under the action by $Spin_7$, we have
$$\tau_1^{-1}(u)=G^{-1}(\tau_1^{-1}(t)).$$

{\bf Corollary 2.4} \ There is no almost  complex structure on the sphere $S^4$.

{\bf Proof}  \ If there is an almost complex structure on the
sphere $S^4$, we have a section $f\colon\; S^4\to CP^3$. As we
know the coholomogy $H^*(CP^3, Z)$ is generated by an element
$a\in H^2(CP^3,Z), \ f^*a=0$, then $f^*(a\cup a)=0$. Let $[\xi]\in
H^4(S^4)$ be a generator, $\tau_1^*[\xi]=\lambda \, a\cup a$.  We
have $f^*\tau_1^*[\xi] =\lambda f^* (a\cup a)=0$, this contradict
to the fact $\tau_1 \circ f =id$. These shows there is no almost
complex structure on the sphere $S^4$.

For Corollary 2.4, see also [2].

By Theorem 2.1, 2.2, 2.3 and  [4] p.342, we have
$$G(2,8)=SO(8)/SO(2)\times SO(6)\approx SO(8)/U(4)\approx SO(7)/U(3),$$
$$CP^3=U(4)/U(1)\times U(3)\approx SO(6)/U(3)\approx SO(5)/U(2),$$
$$CP^1=U(2)/U(1)\times U(1)\approx SO(4)/U(2).$$

\vskip 1cm \centerline{\bf \S 3. The complex structures on $S^6$}
\vskip 0.3cm

In the following we show that  there is no complex structure on
the sphere $S^6$. First we study the differential geometry on the
fibre bundle $\tau \colon\; g(2,8)\to S^6$.

Let $e_1,  e_2,  \cdots,  e_8$ be  orthonormal frame fields on
$R^8$ such that  $e_1 \wedge e_2$ generate a neighborhood of $x$
in $G(2,  8)$.  By $$d(e_1 \wedge e_2)=\sum\limits_{\alpha=3}^8
\omega_1^\alpha E_{i\alpha}+ \sum\limits_{\alpha=3}^8
\omega_1^\alpha E_{2\alpha}, \  \ \omega_i^\alpha =\langle
de_i,e_\alpha\rangle,$$ we know that the elements
$E_{1\alpha}=e_\alpha e_2, \ E_{2\alpha}=e_1e_\alpha, \ \alpha=3,
\cdots, 8, $ can be looked as a basis of $T_{e_1e_2}G(2, 8)$,
$\omega_i^\alpha $ be their dual 1-form. Define the metric on
$G(2, 8)$ by
$$ds^2=2\sum\limits_{i=1}^2 \sum\limits_{\alpha=3}^8 \
(\omega_i^\alpha)^2.$$  Differential $E_{i\alpha}$ we get the
Riemannian connection $\nabla^*$ on $G(2, 8)$,
$$\nabla^* E_{i\alpha}=\sum\limits_{j=1}^2 \ \omega_i^{j}E_{j\alpha}
 +\sum\limits_{\beta=3}^8 \ \omega_\alpha^\beta E_{i\beta}.  $$

Acting $e_1e_2$  on $T_{e_1e_2}G(2,8)$ by Clifford product defines an almost complex $$\widetilde J\colon\;
TG(2,8)\to TG(2,8),$$
$$\widetilde J(E_{1\alpha})=e_1e_2e_\alpha e_2=E_{2\alpha}, \ \
\widetilde J(E_{2\alpha})=e_1e_2e_1e_\alpha =-E_{1\alpha}.$$ It is
easy to see that $(\nabla^*\widetilde J
)E_{i\alpha}=\nabla^*(\widetilde JE_{i\alpha})- \widetilde
J\nabla^* E_{i\alpha}=0$. This shows

{\bf Proposition 3.1} \ The almost complex structure $\widetilde
J$ is integrable and makes $G(2,8)$   a Kaehler manifold.

As is well-known, the Grassmann manifold $G(2,n)$ can be looked as
a complex submanifold of $CP^{n-1}$.

Any $v\in S^6$ defines a complex structure $J_v$ on $R^8$,\
$J_v(\bar e_1) = v, J_v(v)=-\bar e_1$, for any $w\perp \bar e_1,
v, \ J_v(w)$ is defined by  $J_v(w)A = -\bar e_1 vwA$ or
equivalently $ \bar e_1 v\bar e_1wA=\bar e_1J_v(w)A$. For any
$v\in S^6, \ \tau^{-1}(v)=\{uJ_vu \ | \ u\in S^7\}$.

{\bf Proposition 3.2} \ The map $\tau \colon\; G(2,8)\to S^6$ is a
Riemannian  submersion.

{\bf Proof} \  It is easy to see that the vertical subspace
$T_{e_1e_2}\tau ^{-1}(v)$ of the fibre bundle $\tau $  is
generated by
$$e_\alpha J_ve_1 + e_1J_ve_\alpha, \ \ \alpha=3,\cdots,8,$$
where $e_1,e_2=J_ve_1, e_3,\cdots,e_8$ be orthonormal frame fields
along $\tau^{-1}(v)$. The orthogonal subspace of $T_{e_1e_2}\tau
^{-1}(v)$ in $T_{e_1e_2}G(2,8)$ is generated by
$$e_\alpha J_ve_1 -
e_1J_ve_\alpha, \ \alpha=3,\cdots,8.$$ Denote this space by
$T_{e_1e_2}^HG(2,8)$. By
$$(e_\alpha J_ve_1 + e_1J_ve_\alpha)A=0,$$ we have
$$\frac 12(e_\alpha J_ve_1 -
e_1J_ve_\alpha) A= e_\alpha J_ve_1 A=\bar e_1 \tau_*(e_\alpha
J_ve_1)A.$$ The norm of $\frac 12\tau_*(e_\alpha J_ve_1 -
e_1J_ve_\alpha)= \tau_*(e_\alpha J_ve_1)\in T_vS^6$ is $1$. For
$\alpha\neq\beta, \ \tau_*(e_\alpha J_ve_1)\perp \tau_*(e_\beta
J_ve_1)$, we have proved that the map $\tau \colon\; G(2,8)\to
S^6$ is a Riemannian submersion.

Let $TG(2,8)=T^HG(2,8)\oplus T^VG(2,8)$ be the decomposition of
tangent space of $G(2,8), \ T^HG(2,8),T^VG(2,8)$ be horizontal and
vertical spaces respectively. The map $\tau_*:  \
T^H_{e_1e_2}G(2,8)\to T_{\tau(e_1e_2)}S^6$ is an isometric.

{\bf Proposition 3.3} \ (1) The following diagram of maps is
commutative,  $$ \begin{array}{cccccc} T_{e_1e_2}G(2,
8) & \stackrel{\widetilde J} \longrightarrow  & T_{e_1e_2}G(2,   8) \\
\tau_{*} \downarrow \quad & {} & \quad \downarrow \tau_{*} \\
T_vS^6 & \stackrel{J_{e_1e_2}} \longrightarrow & T_vS^6;
\end{array}$$

(2) $\widetilde J\colon\; T^HG(2,8)\to T^HG(2,8),\  \widetilde J:
\ T^VG(2,8)\to T^VG(2,8)$.

{\bf Proof} \ For any $\widetilde X\in T_{e_1e_2}G(2,   8)\subset
\bigwedge^2(R^8), \ \tau_{*}\widetilde X=X$ is defined by
$$\widetilde XA=\bar e_1 XA.$$
(1) follows from the definitions of $\widetilde J$ and $
J_{e_1e_2}$. The proof of (2) is easy.

Any section $f\colon\; S^6\to G(2,8)$ of the fibre bundle $\tau $ defines an almost complex structure $J_f$ on
$S^6, \ J_{f(v)}: \ T_vS^6\to T_vS^6$. As [4], we give

{\bf Definition} \ The section $f$ is holomorphic if
$f_*J_f=\widetilde J f_*$.

From Proposition 3.3 and $\tau f=id$, we have

{\bf Proposition 3.4} \ $J_{f(v)}=\tau_*\widetilde Jf_*.$

For any $X\in T_v^6(S^6), \ f_*(X)=Z_1+Z_2, \ Z_1\in
T_{f(v)}^HG(2,8), Z_2\in  T_{f(v)}^VG(2,8)$, we have
$J_{f(v)}(X)=\tau_*\widetilde J(Z_1).$ On the other hand, the
tangent vector $X$ can be left to a horizontal vector $X^*\in
T_{f(v)}^HG(2,8)$, it is easy to see that $X^*=Z_1$. The almost
complex structure $J_{f(v)}\colon\; T_vS^6\to T_vS^6$ is
determined by the complex structure $\widetilde J$ on
$T^H_{f(v)}G(2,8)$ and the isomorphism $\tau_*\colon\;
T^H_{f(v)}G(2,8)\to T_vS^6$.

By almost complex structure $J_f$, we have $TS^6\otimes
C=T^{(1,0)}S^6\oplus T^{(0,1)}S^6$, where $T^{(1,0)}S^6=\{ X-\sqrt
{-1}J_fX \ | \ X\in TS^6\}, \ T^{(0,1)}S^6=\overline
{T^{(1,0)}S^6}$. The almost complex structure $J_f$ is integrable
if only if
$$[X,Y]\in \Gamma(T^{(1,0)}S^6) \ \mbox{for any} \ X,Y\in
\Gamma(T^{(1,0)}S^6).$$ Similarly, we have decomposition
$TG(2,8)\otimes C = T^{(1,0)}G(2,8)\oplus T^{(0,1)}G(2,8)$ with
respect to the complex structure $\widetilde J$.

 Let $\nabla$ be Riemannian connection on
the sphere $S^6$ and
$$ \$(S^6)=\bigcup_{v\in S^6} \ \{\bar e_1XA \ | \ X\in T_vS^6\otimes
C\}$$ be a vector bundle over $S^6$ which is isomorphic to the
tangent bundle $TS^6\otimes C$. The connection $\nabla$ can be
generalized to the bundle $\$(S^6)$. We have represented the
tangent space of $G(2,8)$ as subspace in $\bigwedge (R^8)$ or
$C\ell_8$. The section $f\colon\; S^6\to G(2,8)$ can  be viewed as
a map $f: \ S^6 \to C\ell_8$ and
$$f_*X= Xf$$
for all $X\in \Gamma (TS^6)$ or $X\in \Gamma (TS^6\otimes C)$.
Note that $\tau_*(Xf) = X$, this can also be written as
$$(Xf)\cdot A= X(fA)=X(\bar e_1vA)=\bar e_1XA.$$

The following is the main result of this section.

{\bf Theorem 3.5} \ The almost complex structure $J_f$ is integrable if and only if the map $f$ is holomorphic.
Then there is no integrable orthogonal complex structure on the sphere $S^6$.

{\bf Proof} \   First assuming  the map $f$ is holomorphic. For
any $X,Y\in \Gamma(T^{(1,0)}S^6)$, we have
$$\widetilde Jf_*(X)=f_*J_{f(v)}X=\sqrt {-1}f_*(X),\ \  \widetilde Jf_*(Y)=\sqrt {-1}f_*(Y), $$
then $ f_*(X),f_*(Y)\in \Gamma(T^{(1,0)}G(2,8))$. The Grassmann
manifold $G(2,8)$ is Kaehlerian, $f_*[X,Y]=[f_*(X),f_*(Y)]\in
\Gamma(T^{(1,0)}G(2,8))|_{f(S^6)}.$  These shows
$$f_*J_f[X,Y]=\widetilde Jf_*[X,Y]=\sqrt {-1}f_*[X,Y].$$
$f\colon\; S^6\to G(2,8)$ is an imbedding, we have $J_f[X,Y]=\sqrt
{-1}[X,Y]$ and the almost complex structure $J_f$ is integrable.

Secondly, assuming the almost complex structure $J_f$ is
integrable. For any $X,Y\in\Gamma(T^{(1,0)}S^6)$, we have
$[X,Y]=\nabla_XY-\nabla_YX\in \Gamma(T^{(1,0)}S^6)$, then
$$(1+\sqrt {-1}f)\bar e_1\nabla_XYA= (1+\sqrt {-1}f)\bar e_1\nabla_YXA.$$
From $f\bar e_1XA=\sqrt {-1}\bar e_1XA$ we have
$$\nabla_Y(f\bar e_1XA)= pr(Yf\cdot \bar e_1YA)+ f\bar
e_1\nabla_YXA=\sqrt {-1}\bar e_1\nabla_YXA,$$ where $pr\colon\;
S^6\times (C\ell^{even}_8\otimes C) \to \$(S^6)$ is a projection
defined naturally. Combine
$$ pr(Yf\cdot \bar
e_1XA)= (\sqrt {-1}-f)\bar e_1\nabla_YXA$$  with $(1+\sqrt
{-1}f)\bar e_1\nabla_XYA= (1+\sqrt {-1}f)\bar e_1\nabla_YXA$, we
have
$$ pr(Yf\cdot \bar
e_1XA)= pr(Xf\cdot \bar e_1YA)= (\sqrt {-1}-f)\bar
e_1\nabla_YXA=(\sqrt {-1}-f)\bar e_1\nabla_YXA.$$ Denote
$\beta(X,Y)= pr(Yf\cdot \bar e_1XA)$.

 $\beta(X,Y)$ is a symmetric
$C^\infty(S^6)$-bilinear form on $\Gamma(T^{(1,0)}S^6)$. We shall
show that $\beta(X,Y)= 0$ is equivalent to the fact that $f$ is
holomorphic. To prove $\beta(X,Y)= 0$ we need only to show
$\beta(X,X)= 0$.

Let $\varepsilon_i=\frac 12(e_{2i-1}-\sqrt{-1}e_{2i}), \ i=1,2,3,
\ e_{2i}= J_fe_{2i-1}, \ \varepsilon_1,
\varepsilon_2,\varepsilon_3$ be local hermitian frame fields for
$T^{(1,0)}S^6$, \
$\varepsilon_i\varepsilon_j=-\varepsilon_j\varepsilon_i, \
\varepsilon_i\overline\varepsilon_i\varepsilon_i =-\varepsilon_i$.
The spinor $\sigma=\varepsilon_1\varepsilon_2\varepsilon_3$ (or
$\varepsilon_1\varepsilon_2\varepsilon_3\overline\varepsilon_1\overline\varepsilon_2\overline\varepsilon_3$)
generates a  spinor bundle on $S^6$ locally.  In the following we
prove
$$\varepsilon_j\nabla_{\varepsilon_j} \sigma =0, \ \ j=1,2,3.$$
The proof is due to  [4], p.335-341, we write here for complete.

Let $\omega_c= (\sqrt{-1})^3e_1e_2\cdots e_6$ be complex volume
element on $TS^6, \ \omega_c\cdot\sigma = -\sigma$. For any $i$,
$$\varepsilon_i\varepsilon_j\nabla_{\varepsilon_j}\sigma=-
\varepsilon_j\varepsilon_i\nabla_{\varepsilon_j}\sigma=-
\varepsilon_j\varepsilon_j\nabla_{\varepsilon_i}\sigma=0,$$ then
 there is a local function $\lambda_{\varepsilon_j}$ such that
$$\varepsilon_j\nabla_{\varepsilon_j}\sigma=
\lambda_{\varepsilon_j}\sigma.$$ On the other hand,
$\nabla_{\varepsilon_j}(\omega_c\sigma)=\omega_c\nabla_{\varepsilon_j}\sigma$,
$$\omega_c\cdot\varepsilon_j\nabla_{\varepsilon_j}\sigma=
(\omega_c\varepsilon_j\omega_c)\nabla_{\varepsilon_j}(\omega_c\sigma)=\varepsilon_j\nabla_{\varepsilon_j}\sigma.$$
These shows $\varepsilon_j\nabla_{\varepsilon_j}\sigma=0$ for any
$j$.

From
$$0=\nabla_{\varepsilon_j}(\varepsilon_j
\sigma)=(\nabla_{\varepsilon_j}\varepsilon_j) \sigma
+\varepsilon_j\nabla_{\varepsilon_j}
\sigma=(\nabla_{\varepsilon_j}\varepsilon_j) \sigma$$ we know that
$\nabla_{\varepsilon_j}\varepsilon_j$ is a local $(1,0)$-frame
field. These shows $\beta(\varepsilon_j,\varepsilon_j)=0$, hence
$\beta(X,Y)=pr(Yf\cdot \bar e_1XA)=0$ for any $X,Y\in
\Gamma(T^{1,0)}S^6)$.

Let $X=(1-\sqrt {-1}J_f)X_1, \ X_1\in \Gamma(TS^6\otimes C), \ 2Yf
= (1-\sqrt {-1}f)Yf+(1+\sqrt {-1}f)Yf$. By $Yf\cdot f =-f\cdot
Yf$, we have
$$(1-\sqrt {-1}f)Yf(1-\sqrt {-1}f)\bar e_1X_1A=(1-\sqrt {-1}f)(1+\sqrt
{-1}f)Yf\bar e_1X_1A=0.$$ Then
$$2pr(Yf\cdot \bar e_1XA)=pr(1+\sqrt{-1}f)^2Yf\bar e_1X_1A=2pr(1+\sqrt{-1}f)Yf\bar e_1X_1A=0.$$

Let $f_*Y=Yf=Z_1+Z_2$ where $Z_1\in \Gamma(T^{(1,0)}G(2,8))$  and
$Z_2\in \Gamma(T^{(0,1)}G(2,8))$. By Proposition 3.3, $Z_2$ is
vertical and we can set $Z_2=(1+\sqrt{-1}f)(X_2J_ve_1+e_1J_vX_2)$
where $X_2\in \Gamma(TS^6\otimes C), \ f(v)=e_1J_ve_1.$ Then for
any $X_1\in \Gamma(TS^6\otimes C)$,
$$pr(1+\sqrt{-1}f)Yf\bar e_1X_1A=pr(1+\sqrt{-1}f)^2(X_2J_ve_1+e_1J_vX_2)\bar e_1X_1A=0.$$
Hence $pr(Yf\cdot \bar e_1XA)=0$  is equivalent to
$$pr(X_2J_ve_1+e_1J_vX_2)\bar e_1X_1A=0.$$
Since $X_1\in \Gamma(TS^6\otimes C)$ is arbitrarily, $X_2=0$, then
$f_*Y$ is a $(1,0)$-form. We have proved that if the almost
complex structure $J_f$ is integrable,
$$\widetilde Jf_*(Y)=\sqrt
{-1}f_*(Y)=f_*J_{f}Y$$ holds for any $Y\in\Gamma(T^{(1,0)}S^6)$.
Then $f$ is holomorphic.

If there is an integrable complex structure on the sphere $S^6$ we
have a holomorphic section $f: \ S^6\to G(2,8)$. Then $S^6$ is
also a Kaehler manifold, this contradict to the fact $H^2(S^6)=0$.

 \vskip 1cm
\centerline{\large \bf References}

\vskip 0.3cm

{\small

\noindent [1] \  M. F. Atiyah, {\it Geometry of Yang-Mills
fields}, Fermi lectures. Scuola Normale, Pisa, 1979.

\noindent [2] \  A. Borel  and  J. P. Serre, {\it Groupes de Lie
et puissances r\'{e}duites de Steenrod}, Amer. J. Math. {\bf
75}(1953), 409-448.

\noindent [3] \ F. R. Harvey, {\it Spinors and calibrations},
Perspectives in Math., 9. Academic Press, New York, 1990.

\noindent [4] \ H. B. Lawson Jr and M. Michelsohn, {\it Spin
geometry}, Princeton University Press, Princeton New Jersey, 1989.

\noindent [5] \ J. W. Zhou and H. Huang, {\it Geometry on
Grassmann Manifolds $G(2, 8)$ and $G(3, 8)$}, Math. J. Okayama
Univ. {\bf 44}(2002), 171-179.

\noindent [6] \  J. W. Zhou, {\it Irreducible Clifford Modules},
Tsukuba J. Math., {\bf 27}(2003), 57-75.

\noindent [7] \ J. W. Zhou, {\it Spinors, Calibrations and
Grassmannians}, Tsukuba Journal of Mathematics, {\bf 27}(2003),
77-97.

\vskip 0.5cm

 \noindent E-mail: \ jwzhou@suda.edu.cn

\end{document}